\documentclass[11pt]{article}
\usepackage{amsfonts}
\usepackage{amssymb}
\usepackage{graphicx}
\usepackage{epsfig}
\usepackage{amsmath}
\usepackage[algo2e,boxed]{algorithm2e}

\newcommand{\ignore}[1]{}
\ignore{The text in between these parentheses is ignored by LaTeX
ignore}

\def\@begintheorem#1#2{\par\bgroup{\sc #1\ #2. }\it\ignorespaces}
\def\@opargbegintheorem#1#2#3{\par\bgroup{\sc #1\ #2\ (#3). } \it\ignorespaces}
\def\@endtheorem{\egroup}
\newtheorem{theorem}{Theorem}[section]
\newtheorem{corollary}[theorem]{Corollary}
\newtheorem{lemma}[theorem]{Lemma}
\newtheorem{proposition}[theorem]{Proposition}
\newtheorem{problem}[theorem]{Problem}
\newtheorem{example}[theorem]{Example}
\newtheorem{algorithm}[theorem]{Algorithm}
\newtheorem{definition}[theorem]{Definition}
\newcommand{\bt}[1]{\begin{theorem}\label{#1}}
\newcommand{\bc}[1]{\begin{corollary}\label{#1}}
\newcommand{\bl}[1]{\begin{lemma}\label{#1}}
\newcommand{\bp}[1]{\begin{proposition}\label{#1}}
\newcommand{\bpro}[1]{\begin{problem}\label{#1}}
\newcommand{\be}[1]{\begin{example}\rm\label{#1}}
\newcommand{\ba}[1]{\begin{algorithm}\rm\label{#1}}
\newcommand{\bd}[1]{\begin{definition}\rm\label{#1}}
\newcommand{\bpr}{\par{\it Proof}. \ignorespaces}
\newcommand{\et}{\end{theorem}}
\newcommand{\ec}{\end{corollary}}
\newcommand{\el}{\end{lemma}}
\newcommand{\ep}{\end{proposition}}
\newcommand{\epro}{\end{problem}}
\newcommand{\ee}{\end{example}}
\newcommand{\ea}{\end{algorithm}}
\newcommand{\ed}{\end{definition}}
\newcommand{\epr}{{\ \vbox{\hrule\hbox{%
\vrule height1.3ex\hskip0.8ex\vrule}\hrule}}\\\par}

\def\N{\mathbb{N}}
\def\R{\mathbb{R}}
\def\Z{\mathbb{Z}}
\def\F{\mathbb{F}}

\def \diag {{\rm diag}}
\def \rank {{\rm rank}}
\def \supp {{\rm supp}}

\def \l {\langle}
\def \r {\rangle}

\def \B {{{\mathcal B}}}
\def \I {{{\mathcal I}}}
\def \E {{{\mathcal E}}}

\begin{document}

\title{\bf Nonlinear Matroid Optimization
\break and Experimental Design}
\author{\qquad\qquad
Yael Berstein
\qquad
Jon Lee
\qquad
Hugo Maruri-Aguilar
\qquad
Shmuel Onn
\qquad\qquad
\and
\qquad\qquad
Eva Riccomagno
\qquad
Robert Weismantel
\qquad
Henry Wynn
\qquad\qquad
}
\date{}
\maketitle

\begin{abstract}
We study the problem of optimizing nonlinear objective functions
over matroids presented by oracles or explicitly. Such functions
can be interpreted as the balancing of multi-criteria optimization.
We provide a combinatorial polynomial time algorithm for arbitrary
oracle-presented matroids, that makes repeated use of matroid
intersection, and an algebraic algorithm for vectorial matroids.

Our work is partly motivated by applications to
minimum-aberration model-fitting in experimental design
in statistics, which we discuss and demonstrate in detail.
\end{abstract}

\section{Introduction}

In this article, partly motivated by applications to minimum-aberration
model-fitting in experimental design, which will be discussed briefly
at the end of this introduction and in detail in \S5,
we study the problem of optimizing an arbitrary nonlinear function
over a matroid, set as follows.

\vskip.2cm\noindent {\bf Nonlinear Matroid Optimization}. Given
matroid $M$ on ground set $N:=\{1,\dots,n\}$, integer weight vectors
$w_i=(w_{i,1},\dots,w_{i,n})\in\Z^n$
for $i=1,\dots,d$, and function $f:\R^d\rightarrow\R$,
find a matroid base $B\in\B(M)\subset 2^N$ minimizing the ``balancing"
by $f$ of the $d$ weights $w_i(B):=\sum_{j\in B} w_{i,j}$ of base $B$,
$$f(w_1(B),\dots,w_d(B))\ =\
f\left(\sum_{j\in B}w_{1,j},\dots,\sum_{j\in B}w_{d,j}\right)\ .$$

\vskip.2cm\noindent
(All necessary basics about matroid theory are provided in \S\ref{Preliminaries}.
For more details consult \cite{Oxley} or \cite{Wel}.)

Nonlinear matroid optimization can be interpreted as
{\em multi-criteria matroid optimization}: the $d$ given
weight vectors $w_1,\dots,w_d$ represent $d$ different criteria,
where the value of base $B\in\B(M)$ under criterion $i$ is its
$i$-th weight $w_i(B)=\sum_{j\in B}w_{i,j}$, and where the objective
is to minimize the ``balancing" $f(w_1(B),\dots,w_d(B))$ of the $d$
given criteria by the given function $f$.

In fact, we have a hierarchy of problems of increasing generality,
parameterized by the number $d$ of weight vectors.
At the bottom lies standard linear matroid optimization,
recovered with $d=1$ and $f$ the identity on $\R$.
At the top lies the problem of maximizing an arbitrary function over
the set of bases, with $d=n$ and $w_i={\bf 1}_i$ the $i$-th standard unit
vector in $\R^n$ for all $i$, see Proposition \ref{Intractability2} below.

It will be often convenient to collect the weight vectors in a $d\times n$
matrix $W$. Thus, the $i$-th row of this matrix is the $i$-th weight vector
$w_i=(w_{i,1},\dots,w_{i,n})$. For each subset $S\subseteq N$ we define its
$W$-\emph{profile} to be
\[
W(S)\ :=\ (w_1(S),\dots,w_d(S))
\ =\ \left(\sum_{j\in S}w_{1,j},\dots,\sum_{j\in S}w_{d,j}\right)\ \in\ \Z^d~,
\]
which is the vector giving the value of $S$ under each of the $d$ weight vectors.
The nonlinear matroid optimization then asks for a matroid base $B\in\B(M)$
minimizing the objective function $f(W(B))$.

The computational complexity of nonlinear matroid optimization
depends on the number $d$ of weight vectors, on the representation
of weights (binary $\l w_{i,j}\r$ versus unary $|w_{i,j}|$, see
\S\ref{Preliminaries}), on the type of function $f$ and its presentation,
and on the type of matroid and its presentation. We will be able to
handle an arbitrary function $f$ presented by a {\em comparison oracle}
that, queried on $x,y\in\Z^d$, asserts whether or not $f(x)\leq f(y)$,
and an arbitrary matroid presented by an {\em independence oracle}, that,
queried on $I\subseteq N$, asserts whether or not $I$ is {\em independent} in $M$,
that is, whether $I\subseteq B$ for some base $B\in\B(M)$, see \S\ref{Preliminaries}.
These are very broad presentations that reveal little information (per query)
on the function and matroid, making our problem setting rather expressive
but difficult for achieving strong complexity results.

Standard linear matroid optimization is well known to be efficiently
solvable by the greedy algorithm. Nonlinear matroid optimization
with $d$ fixed and $f$ \emph{concave} turns out to be solvable in
polynomial time as well \cite{HT,Onn}. In fact, using sophisticated
geometric methods, the nonlinear optimization problem
with $d$ fixed and $f$ concave has been recently shown to be efficiently
solvable for bipartite matching \cite{BO} and for broader classes of
combinatorial optimization problems \cite{OR}.
Therein, maximization rather than minimization form is used,
and hence convex rather than concave functions are considered.

However, generally, nonlinear matroid optimization is intractable, even
for uniform matroids. In particular, if $d$ is variable then exponential
time is needed even for $\{0,1\}$-valued weights, and if the weights
are presented in binary, then exponential time is needed even for
fixed dimension $d=1$. See Propositions \ref{Intractability1},
\ref{Intractability2} and \ref{Intractability3} in the sequel
for various intractability statements.

In spite of these difficulties, here we establish the
efficient solvability of the problem as follows.

\bt{Arbitrary}
For every fixed $d$ and $p$, there is an algorithm that, given
a matroid $M$ presented by an independence oracle on  the $n$-element ground set $N$,
integers $a_1,\dots,a_p\in\Z$,
weight vectors $w_1,\dots,w_d\in\{a_1,\dots,a_p\}^n$, and function
$f:\R^d\rightarrow\R$ presented by a comparison oracle, solves the nonlinear
matroid optimization problem in time that is polynomial in $n$ and $\max \l a_i\r$.
\et

We also state the following natural immediate corollary
concerning $\{0,1,\dots,p\}$-valued weights.

\bc{Corollary}
For every fixed $d$ and $p$, there is an algorithm that, given $n$-element matroid $M$
presented by an independence oracle, $w_1,\dots,w_d\in\{0,1,\dots, p\}^n$,
and function $f:\R^d\rightarrow\R$ presented by a comparison oracle,
solves the nonlinear matroid optimization problem in time polynomial in $n$.
\ec

The algorithm establishing Theorem \ref{Arbitrary} is combinatorial
and makes repeated use of matroid intersection (see e.g. \cite{Lee}
or \cite{LeeRyan}, and see \cite{HL} for another recent interesting application
of matroid intersection). In fact, it invokes the matroid intersection
algorithm roughly $n^{p^d}$ times, and hence it is quite heavy.
However, most matroids appearing in practice, including graphic matroids
and those arising in the applications to experimental design to be
discussed later, are vectorial. Therefore, we also provide another, more efficient,
linear-algebraic algorithm for vectorial matroids. Moreover, this algorithm
applies to weights with unlimited number (rather than fixed number $p$)
of different values $w_{i,j}$ of entries.

\bt{Vectorial}
For every fixed $d$, there is an algorithm that, given integer $m\times n$
matrix $A$, weight vectors $w_1,\dots,w_d\in\Z^n$, and
function $f:\R^d\rightarrow\R$ presented by a comparison oracle,
solves the nonlinear optimization problem over the vectorial matroid
of $A$ in time polynomial in $\l A\r$ and $\max|w_{i,j}|$.
\et

A specific application that can be solved by either the combinatorial
algorithm underlying Theorem \ref{Arbitrary} or the more efficient
linear-algebraic algorithm underlying Theorem \ref{Vectorial} is the following.

\be{Graphic}{\sc (minimum-norm spanning tree).}
Fix any positive integer $d$. Let $G$ be any connected graph with edge set
$E=\{e_1,\dots,e_n\}$, and let $w_1,\dots,w_d\in\Z^n$ be weight vectors with
$w_{i,j}$ representing the values of edge $e_j$ under the $i$-th criterion.
Let $A$ be the vertex-arc incidence matrix of an arbitrary orientation of $G$.
Then the vectorial matroid of $A$ is the graphic matroid of $G$ whose bases are
the spanning trees of $G$. Now fix also any $q$ that is either a positive integer
or $\infty$, and let $f:\R^d\rightarrow\R$ be the $l_q$~norm given by
$\|x\|_q:=(\sum_{i=1}^d |x_i|^q)^{1\over q}$ for finite $q$ and
$\|x\|_{\infty}:=\max_{i=1}^d |x_i|$. Note that a comparison oracle
for $f=\|x\|_q$ is easily and efficiently realizable.
Then Theorems \ref{Arbitrary} and \ref{Vectorial} assure that a spanning tree
$T$ of $G$ minimizing the $l_q$~norm of the multi-criteria vector, given by
$$\left\|\left(w_1(T),\dots,w_d(T)\right)\right\|_q\ =\
\left\|W(T)\right\|_q
$$
is computable in time polynomial in $n$ and $\max|w_{i,j}|$. Note that if P$\neq$NP
then the problem is {\em not} solvable in time polynomial in the binary length
$\l w_{i,j}\r$ even for the graph obtained from a path by replacing every
edge by two parallel copies, see Proposition \ref{Intractability3} below.
\ee

\vskip.4cm\noindent{\bf Experimental Design.}
We conclude the introduction with a brief discussion of the application
to experimental design, elaborated in detail in \S5. The general framework is as follows.
We are interested in learning an unknown system whose output $y$ is an unknown
function $\Phi$ of a multivariate input $x=(x_1,\dots,x_k)\in\R^k$. It is customary to
call the input variables $x_i$ {\em factors} of the system.
In order to learn the system, we perform several experiments.
Each experiment $i$ is determined by a point $p_i=(p_{i,1},\dots,p_{i,k})$
and consists of feeding the system with input $x:=p_i$ and measuring the
corresponding output $y_i=\Phi(p_i)$. Based on these experiments, we wish to
{\em fit a model} for the system, namely, determine an estimation $\hat \Phi$ of the
function $\Phi$, that satisfies the following properties:
\begin{itemize}
\item
It lies in a prescribed class of functions;
\item
It is consistent with the outcomes of the experiments;
\item
It minimizes the {\em aberration} - a suitable criterion -
among models in the class.
\end{itemize}
More detailed discussion and precise definitions will be given in
\S5. We have the following broad corollary of Theorems \ref{Arbitrary}
and \ref{Vectorial}; see \S5 for the precise statement and its various
practical specializations to concrete aberration criteria useful
in optimal model-fitting in experimental design.

\bc{ExperimentalDesign}
An aberration-minimum multivariate-polynomial model is polynomial time computable.
\ec

The article proceeds as follows. In \S2 we set some notation and preliminaries,
make some preparations for the algorithms in following sections, and demonstrate
various intractability limitations on nonlinear matroid optimization.
In \S3 we discuss arbitrary matroids presented by oracles,
and provide the combinatorial algorithm for nonlinear matroid optimization,
thereby establishing Theorem \ref{Arbitrary}. In \S4 we provide the
more efficient algebraic algorithm for nonlinear optimization over vectorial
matroids, thereby proving Theorem \ref{Vectorial}. We conclude in \S5 with
a detailed discussion of the experimental design application and prove Corollary
\ref{ExperimentalDesign} and its various practical specializations.
Readers interested mostly in experimental design can go directly
to \S5, where the minimum-aberration model-fitting problem is reduced to
nonlinear optimization over a suitable matroid, and where each of the
algorithms developed in \S3-4 can be invoked as a black box.

\section{Preliminaries, Preparation and Limitations}

\subsection{Preliminaries}
\label{Preliminaries}

We use $\R$ for the reals, $\Z$ for the integers, and $\N$ for the nonnegative integers.
The $i$-th standard unit vector in $\R^n$ is denoted by ${\bf 1}_i$. The {\em support}
of $x\in\R^n$ is the index set $\supp(x):=\{j:x_j\neq 0\}$ of nonzero entries of $x$.
The integer lattice $\Z^n$ is naturally embedded in $\R^n$.
Vectors will be interpreted as either row or column vectors
interchangeably -- this will be relevant only when such vectors are
multiplied by matrices -- in which case the correct interpretation
will be clear from the context. We often collect a sequence of
vectors designated by a low case letter as the rows of a matrix
designated by the corresponding upper case letter. Thus,
our weight vectors $w_i=(w_{i,1},\dots,w_{i,n})$, $i=1,\dots,d$
are arranged in a $d\times n$ matrix $W$, and our design points
$p_i=(p_{i,1},\dots,p_{i,k})$, $i=1,\dots,m$ are arranged in
an $m\times k$ matrix $P$. The space $\R^n$ is endowed with the standard
inner product which, for $w,x\in\R^n$, is $w\cdot x:=\sum_{i=1}^n w_ix_i$.
Vectors $w$ in $\R^n$ are also regarded as linear functions
on $\R^n$ via the inner product $w\cdot x$. So we refer to elements
of $\R^n$ as points, vectors, or linear functions, as is
appropriate from the context. We write ${\bf 1}$ for the vector
with all entries equal to $1$, of dimension clear from the context.

Our algorithms are applied to rational data only, and the time complexity
is as in the standard Turing machine model, see e.g. \cite{AHU,GJ,Sch}.
The input typically consists of rational (usually integer) numbers,
vectors, matrices, and finite sets of such objects.
The {\em binary length} of an integer number $z\in\Z$ is defined
to be the number of bits in its binary representation,
$\l z \r:=1+\lceil \log_2(|z|+1)\rceil$ (with the extra bit for the sign).
The length of a rational number presented as a fraction $r={p\over q}$
with $p,q\in\Z$ is $\l r \r:=\l p\r + \l q\r$. The length of an
$m\times n$ matrix $A$ (or a vector) is the sum
$\l A \r:=\sum_{i,j}\l a_{i,j}\r$ of the lengths of its entries.
Note that the length of $A$ is no smaller than the number of entries,
$\l A \r\geq mn$. Therefore, when $A$ is, say, part of an input to an
algorithm, with $m,n$ variable, the length $\l A \r$ already incorporates
$mn$, and so we will typically not account additionally for $m,n$ directly.
But sometimes we will also emphasize $n$ as part of the input length.
Similarly, the length of a finite set $S$ of numbers, vectors or matrices
is the sum of lengths of its elements and hence, since $\l S\r\geq |S|$,
automatically accounts for its cardinality.
Some input numbers affect the running time of some algorithms through their
unary presentation, resulting in so-called ``pseudo-polynomial" running time.
The {\em unary length} of an integer number $z\in\Z$ is the number $|z|+1$ of
bits in its unary representation (again, an extra bit for the sign). The unary
length of a rational number, vector, matrix, or finite set of such objects
are defined again as the sums of lengths of their numerical constituents,
and is again no smaller than the number of such numerical constituents.
Often part of the input is presented by oracles. Then the running time counts
also the number of oracle queries. An oracle algorithm is {\em polynomial time}
if its running time, including the number of oracle queries,
is polynomial in the length of the input.

Next, we make some basic definitions concerning matroids and set some
associated notation. For matroid theory see \cite{Oxley} or \cite{Wel}.
A {\em matroid} $M$ is described by giving a finite \emph{ground set}
$\E(M)$ and a nonempty set $\B(M)$ of subsets  of $\E(M)$ called the set of
{\em bases} of $M$, such that for every $B,B'\in\B(M)$ and $i\in
B\setminus B'$ there is an $i'\in B'$ such that
$B\setminus\{i\}\cup\{i'\}\in\B(M)$. A subset $I$ of $\E(M)$ is
{\em independent} in $M$ if $I\subseteq B$ for some base $B\in\B(M)$.
The set of independent sets of $M$ is denoted by $\I(M)$.
Often it is convenient to let ${\mathcal E}(M)=N=\{1,2,\ldots,n\}$
for some positive integer $n$.
The {\em rank} $\rank(S)$ of a subset $S\subseteq \E(M)$ is the maximum
cardinality $|I|$ of an $I\in\I(M)$ that is contained in $S$. The
{\em rank} of the matroid is $\rank(M):=\rank(\E(M))$ and is equal to the
cardinality $|B|$ of every base $B\in\B(M)$. An {\em independence
oracle} for $M$ is one that, queried on $I\subseteq \E(M)$, asserts
whether $I$ is in $\I(M)$. An independence oracle allows us to
compute the rank of every $S\subseteq \E(M)$ as follows. Start with
$I:=\emptyset$. For each $j\in S$ do (in any order): if $I\cup\{j\}$
is in $\I(M)$ then set $I:=I\setminus\{j\}$. Output
$\rank(S):=|I|$.

An important  example of a matroid is the {\em graphic
matroid} $M(G)$ of a graph $G$ with $\E(M(G))$ being the
edge set of $G$ and $\B(M(G))$
equal to the set of edge sets of spanning forests of $G$.
A further key example is the {\em vectorial matroid}
of an $m\times n$ matrix $A$ (over a field $\F$) with $\B(M)$ the set of subsets of
indices of maximal linearly-independent subsets of columns of $A$.
Throughout, when treating complexity issues, we assume that $A$ has
integer components, and that we do arithmetic over the rationals.
A {\em uniform matroid} is any matroid that is isomorphic to the matroid
${\mathcal U}_{m,n}$ having  ground set $N=\{1,\dots,n\}$ and
having all $m$-subsets of $N$ as bases.
If $M_1$ and $M_2$ are matroids with disjoint ground sets, then their \emph{direct sum}
$M_1\oplus M_2$ has $\E(M_1\oplus M_2):=\E(M_1)\uplus \E(M_2)$ and
$\B(M_1\oplus M_2):=\{B_1\uplus B_2:B_1\in\B(M_1), B_2\in \B(M_2)\}$.
A \emph{partition matroid} is any matroid that is a direct sum of uniform matroids
$\oplus_{i=1}^r {\mathcal U}_{m_i,n_i}$
(with the ground sets of the ${\mathcal U}_{m_i,n_i}$ labeled
to be pairwise disjoint). If $m_i=m$ and all $n_i=k$ for all $i$,
then we write ${\mathcal U}^r_{m,k}$ for this $r$-fold sum of ${\mathcal U}_{m,k}$.

It is well known that all graphic matroids and uniform matroids are vectorial.
Furthermore, any partition matroid that is the direct sum
$\oplus_{i=1}^r {\mathcal U}_{1,n_i}$ of rank-1 uniform matroids
is graphic: it is the matroid of the graph obtained from any forest on $r$
edges by replacing the $i$-th edge by $n_i$ parallel copies for $i=1,\dots,r$.
In particular, ${\mathcal U}^r_{1,k}$ is a graphic matroid for any positive $k$ and $r$.

\subsection{Preparation}
\label{Preparation}

Here we provide preparatory ingredients which will be used in
algorithms in \S3-4. We begin by showing that finding an
optimal base for a nonlinear matroid optimization problem can be
reduced to finding the optimal objective function {\em value} of a
small number of subproblems. Consider data for a nonlinear matroid
optimization problem, consisting of a matroid $M$, weight
vectors $w_1,\dots,w_d\in\Z^n$, and function
$f:\R^d\rightarrow\R$. Each subset $S\subseteq N$ gives a
\emph{subproblem} of nonlinear matroid optimization as follows. The matroid
of the subproblem is the {\em restriction} of $M$ to $S$, that is,
the matroid $M.S$ on ground set $S$ in which a subset $I\subseteq S$ is
independent if and only if it is independent in $M$. Note that an
independence oracle for the restriction matroid $M.S$ is realizable at once
from that of $M$. The weight vectors of the subproblem are the
restrictions of the original weight vectors to $S$. The function
$f:\R^d\rightarrow\R$ in the subproblem is the same as in the
original problem. We have the following useful statement.

\bl{Objective} The nonlinear matroid optimization problem of finding
an optimal base of an $n$-element matroid is reducible in time
polynomial in $n$ to finding the optimal objective value of $n+1$
subproblems. \el
\bpr Denote by $f^*(S)$ the optimal objective
function value of the subproblem on $S\subseteq N$. Now compute
the optimal objective function value for $n+1$ such subproblems as
follows.

\begin{algorithm2e}[H]
Start with $S:=N$\;
Compute $m:=\rank(S)$\;
Compute the optimal value $f^*:=f^*(N)$ of the original problem\;
\For{j=1,2,\dots,n}{
Compute $\rank(S\setminus\{j\})$\;
Compute the optimal subproblem value $f^*(S\setminus\{j\})$\;
\lIf {$\rank(S\setminus\{j\})=m$ and
$f^*(S\setminus\{j\})=f^*$}{set $S:=S\setminus\{j\}$}\;
}
{\bf return} $B:=S$\;
\end{algorithm2e}

\noindent It is not hard to verify that the set $B$
obtained is indeed an optimal base for the original problem. \epr

We record for later the following statement,
which follows directly from the definitions.
\bp{BaseWeightVector}
Consider $n$-element matroid $M$, weights $w_1,\dots,w_d\in\Z^n$,
and $f:\R^d\rightarrow \R$. Put
\begin{equation}\label{e_U}
U\ :=\ \{W(B)\ :\ B\in\B(M)\}\ =\ \{u\in \Z^d \ :\  u=W(B)\ \mbox{for some
base}\ B\in\B(M)\}\ .
\end{equation}
Then the optimal objective value of the corresponding
nonlinear matroid optimization problem satisfies
$$f^*\ =\ \min_{u\in U} f(u)\ .$$
\ep Thus, the problem of computing $f^*$ reduces to that of
constructing the set $U$ of $W$-profiles of bases. This really is
the crucial component of the solution of the nonlinear matroid
optimization problem. While the cardinality of $U=\{W(B)\ :\ B\in\B(M)\}$
may be polynomial under suitable assumptions on the data, its direct
computation by computing $W(B)$ for every base is prohibitive since
the number of matroid bases is typically exponential in $n$.
Instead, we will construct a finite superset $Z$ of ``potential"
$W$-profiles of bases, satisfying $U\subseteq Z\subset\Z^d$, and
then filter $U$ out of $Z$. However, this is not an easy task
either: deciding if a given $u\in\Z^d$ satisfies $u=W(B)$ for some
base $B$ is NP-complete already for fixed $d=1$ and uniform
matroids or partition matroids that are the direct sums of rank-1 uniform matroids,
see Proposition \ref{Intractability3} below. The way the filtration of $U$ out of
$Z$ is done is precisely the key difference between our two
algorithms for nonlinear matroid optimization in \S3 and \S4.

\subsection{Limitations}
\label{Limitations}

Here we provide various intractability statements about the
nonlinear matroid optimization problem and some of its relatives. In
particular, we show that, if the weights are encoded in binary, then
even the $1$-dimensional problem, namely with fixed $d=1$, over any
matroid given explicitly or by an oracle, requires to examine the
objective value of {\em every} base, and hence cannot be solved in
polynomial time.

It is convenient to define a certain class of rank $k$ matroids on $2k$ elements,
so as to generalize the uniform matroid ${\mathcal U}_{k,2k}$
and the partition matroids ${\mathcal U}^k_{1,2}$.
For some $r$, $1\le r\le k$,
partition $K=\{1,\ldots,k\}$ into $r$ parts, as $\uplus_{i=1}^r K_i$.
Correspondingly, let $\bar{K}_i := \{ \bar{j}\ :\ j\in K_i \}$.
Let $k_i=|K_i|=|\bar{K}_i|$. For $i=1,\ldots,r$, let ${\mathcal M}_i$ be a
uniform matroid of rank $k_i$ on ground set $K_i \cup \bar{K}_i$
(so ${\mathcal M}_i\cong {\mathcal U}_{k_i,2k_i}$).
Let ${\mathcal M}_{r,k} := \oplus_{i=1}^r {\mathcal M}_i$~.
Observe that ${\mathcal M}_{1,k} \cong {\mathcal U}_{k,2k}$ (which is uniform)
and ${\mathcal M}_{k,k} \cong {\mathcal U}^k_{1,2}$ (which is graphic),
so \emph{hardness results with respect to the matroid classes ${\mathcal M}_{r,k}$
apply to uniform and graphic matroids.}
We note and will soon use that for any $r$, $1\le r\le k$, $|\B({\mathcal M}_{r,k})|$
is not bounded by any polynomial in $|\E({\mathcal M}_{r,k})|=2k$.

\bp{Intractability1} Computing the optimal objective value of the
$1$-dimensional nonlinear matroid optimization problem
over matroids $M$, given explicitly or by an oracle, on $n$-element ground sets,
 and a univariate function
$f$ presented by a comparison oracle, cannot be done in polynomial time.
In particular, with the single
weight vector $w:=(1,2,4,\dots,2^{n-1})$,
solution of the nonlinear matroid optimization problem
requires examining $f(W(B))$ for each of the $|\B(M)|$ bases of $M$.
In particular, for each $r$, $1\le r\le k$, the problem cannot be solved
in polynomial time for the class of matroids ${\mathcal M}_{r,k}$~.
\ep

\bpr
The weights $w(S)=\sum_{j\in S}2^{j-1}$ of
the $2^n$ subsets $S\subseteq N$ attain precisely all $2^n$ distinct
values $0,1,\dots,2^n-1$. Since the function $f$ is arbitrary, this
implies that the objective value $f(W(B))$ of each base $B$ can be
arbitrary. Therefore, if the value $f(W(B))$ of some base $B$ is
not examined, it may be that this value is the unique minimum one
and the nonlinear matroid optimization problem cannot be correctly
solved. The final remark of the proposition follows since, for every $r$,
$|\B({\mathcal M}_{r,k})|$ is not bounded
by any polynomial in $|\E({\mathcal M}_{r,k})|=2k$.
\epr

\bp{Intractability2}
Computing the optimal objective value of the
nonlinear matroid optimization problem in variable dimension $d=n$,
over any matroid $M$, with $\{0,1\}$-valued weights, the $i$-th weight being
the standard unit vector $w_i:={\bf 1_i}$ in $\R^n$ for all $i$, and with
$f:\R^n\rightarrow \R$ a function presented by a comparison
oracle, requires  examining $f(W(B))$ for each of the $|\B(M)|$ bases of $M$.
In particular, for each $r$, $1\le r\le k$, the problem cannot be solved
in polynomial time for the class of matroids ${\mathcal M}_{r,k}$~.
\ep

\bpr
The $W$-profiles $W(B)=(w_1(B),\dots,w_n(B))$ of the $2^n$ subsets
$B\subseteq N$ attain precisely all $2^n$ distinct vectors in
$\{0,1\}^n$. Since the function $f$ is arbitrary, this implies that
the objective value $f(W(B))$ of each base can be arbitrary. The
rest of the argument is as in the proof of Proposition
\ref{Intractability1}. \epr

Consider nonlinear matroid optimization with a matroid $M$,
weights $w_1,\dots,w_d\in\Z^n$, and function $f:\R^d\rightarrow \R$.
As explained in \S \ref{Preparation}, a crucial component in
solving the problem is to identify $W$-profiles of bases; that is,
points $u\in\Z^d$ satisfying $u=W(B)$ for some $B\in\B(M)$. The
following proposition shows that, with binary encoded weights, both
the nonlinear matroid optimization problem and the problem of
deciding if a given $u$ is a $W$-profile of some base, are
intractable already for fixed $d=1$ and uniform matroids or
partition matroids that are the direct
sums of rank-1 uniform matroids.

\bp{Intractability3}
Given matroid $M$, a single nonnegative weight vector $w\in\N^n$,
and nonnegative integer $u\in\N$, encoded in binary, the following problems are
NP-complete, already when restricted to the class of matroids
${\mathcal M}_{r,k}$, for any $r$, $1\le r \le k$:
\begin{enumerate}
\item
Determining whether $u=W(B)=\sum_{j\in B}w_j$ for some base $B\in\B(M)$.
\item
Determining whether the optimal objective value is zero for the $1$-dimensional
nonlinear matroid optimization problem over $M$, with the explicit convex
univariate function $f(y):=(y-u)^2$.
\end{enumerate}
\ep

\bpr
The NP-complete {\em subset-sum problem} is to decide,
given $a_0,a_1,\dots,a_k\in\N$, whether there is a subset
$S\subseteq K=\{1,\dots,k\}$ with $\sum_{j\in S}a_j=a_0$.
Given such
$a_i$, let $u:=a_0$, and let
$w:=(a_1,\dots,a_k,0,\dots,0)\in \N^K\times\N^{\bar{K}}$.
Then, for any $r$, a base $B$ of ${\mathcal M}_{r,k}$
satisfies $W(B)=u$ if and only
if $S:=B\cap K$ satisfies $\sum_{j\in S}a_j=a_0$.
This
reduces the subset-sum problem to the problem considered in the
first part of the proposition, showing that it is indeed
NP-complete.

For the second part, note that the objective value
$f(W(B))=(W(B)-u)^2$ of every base $B$ is nonnegative, and $B$ has
$f(W(B))=0$ if and only if $W(B)=u$. Thus, the optimal objective
value is zero if and only if there is a base with $W(B)=u$. So the
problem in the first part of the proposition reduces to the problem
in the second part, showing the latter to be NP-complete as well.
\epr

\section{Arbitrary Matroids}

In this section we develop a combinatorial algorithm for nonlinear matroid optimization,
that runs in polynomial time for any matroid presented by an independence oracle,
provided that the number $p$ of distinct values taken by the entries $w_{i,j}$
of the weight vectors is fixed. In particular, the algorithm applies to
$\{0,1\}$-valued weight vectors as well as to $\{0,1,\dots,p\}$-valued
weight vectors for any fixed $p$.

As explained in \S \ref{Preparation}, we will filter the set
$U=\{W(B)\ :\ B\in\B(M)\}$ of $W$-profiles of bases out of a suitable superset $Z$.
For this, we next show how to efficiently decide if a given $u\in\Z^d$ satisfies
$u=W(B)$ for some $B\in\B(M)$. We start with  $\{0,1\}$-valued
$W$-profiles with pairwise disjoint supports.

\bl{Disjoint}\label{Disjoint} There is an algorithm that, given matroid
$M$ presented by an independence oracle on the $n$-element ground set $N$,
weight vectors $w_1,\dots,w_d\in\{0,1\}^n$ with pairwise disjoint
supports, and $u\in\N^d$, determines if $M$ has a base $B$ with
$W$-profile $W(B)=u$, in time polynomial in $n$ and
$\l u\r$. \el

\bpr
For each base $B$ of $M$ and for each $i=1,\dots,d$
we have $w_i(B)=|B\cap\supp(w_i)|$. Therefore, a base $B$ has
$W(B)=u$ if and only if $|B\cap\supp(w_i)|= u_i$ for $i=1,\dots,d$.
So we may and do assume $\sum_{i=1}^d u_i\leq \rank(M)$ and
$u_i \leq |\supp(w_i)|$ for all $i$ else $M$ has no base with $W(B)=u$. Let
\begin{align*}
\B' :=& \biggl\{B\subseteq N\ :\ |B\cap\supp(w_i)|=u_i~,\quad i=1,\dots,d~,\\
& \textstyle \quad\quad\quad\quad\quad\,
\left|B\cap \left(N\setminus \bigcup_{i=1}^d\supp(w_i)\right)\right|
=\rank(M)-\sum_{i=1}^d u_i\biggr\}~.
\end{align*}
It is easy to see that $\B'=\B(M')$ is the set of bases of a partition matroid $M'$ on
ground set $N$, for which an independence oracle is efficiently realizable.
Moreover, $M$ has a base $B$ with $W(B)=u$ if and only if $\B(M)\cap\B(M')$
is nonempty. These observations justify the following algorithm:

\begin{algorithm2e}[H]
\lIf {$\sum_{i=1}^d u_i > \rank(M)$ or
$u_i > |\supp(w_i)|$ for some $i=1,\ldots,d$}{{\bf return} NO\;}
Determine if  $\B(M)\cap\B(M')$ is nonempty by
computing a max-cardinality $S\in \I(M) \cap \I(M')$\;
\eIf {$|S|= \rank(M)$\;}{{\bf return} YES and $B:=S$\;}{{\bf return} NO\;}
\end{algorithm2e}

Computing a max-cardinality $S\in \I(M) \cap \I(M')$ can be  efficiently carried out
using a maximum-cardinality matroid-intersection algorithm, see e.g.
\cite{Lee} or \cite{LeeRyan} and references therein. \epr

Next we consider weight vectors for which the number $p$ of distinct $w_{i,j}$ values
is fixed. So, we assume that $w_1,\dots,w_d\in\{a_1,\dots,a_p\}^n$ for arbitrary
given integer numbers $a_1,\dots,a_p$. Note that the $a_i$ can vary and be very
large, since they affect the running time through their binary length $\l a_i\r$.

\bl{Unary} For every fixed $d$ and $p$, there is an algorithm that,
given matroid $M$ presented by an independence oracle on ground set $N$,
integers $a_1,\dots,a_p$, weight vectors
$w_1,\dots,w_d\in\{a_1,\dots,a_p\}^n$, and $u\in\N^d$, decides if
$M$ has a base  $B$ with $W(B)=u$, in time polynomial in $n$, $\max
\l a_i \r$, and $\l u \r$. \el

\bpr
Let $V:=\{a_1,\dots,a_p\}^d$ and let $P$ be the $d\times p^d$ \emph{pattern matrix}
having columns that are all of the $p^d$ points in $V$. Let, as usual,
$W$ be the $d\times n$ matrix with rows $w_1,\dots,w_d$.
For $j=1,\dots,n$, let $w^j$ denote the $j$-th column of $W$.
We will exploit the fact that, no matter how large $n$ is, the columns $w^j$ of $W$ all
lie in the fixed set $V$ -- so the number of distinct columns $w^j$ of $W$ is limited.

Define a $p^d\times n$ \emph{selector matrix} ${\widehat W}$,
having rows ${\widehat w}_v$ indexed by $V$, columns ${\widehat w}^j$ indexed by $N$, and
$${\widehat w}_{v,j}\ :=\
           \left\{
             \begin{array}{ll}
               1, & \hbox{if $w^j=v$~;} \\
               0, & \hbox{otherwise,}
             \end{array}
           \right.
$$
for $v\in V$, $j\in N$. Note that each column $\widehat{w}^j$ of ${\widehat W}$
is a standard unit vector, selecting the unique pattern (i.e., column) of $V$
that agrees with the column $w^j$ of $W$. It should be clear that the rows
of ${\widehat W}$, namely the ${\widehat w}_v$, lie in $\{0,1\}^n$ and have pairwise
disjoint supports -- this will enable us to appeal to Lemma \ref{Disjoint}.
We observe and will make use of the fact that the weight matrix $W$ factors as
$W=P{\widehat W}$. Therefore, the $W$-profile and $\widehat{W}$-profile of a base $B$
satisfy $W(B)= P{\widehat W}(B)$~. This implies that there exists a base $B\in\B(M)$
with $W(B)=u$ if and only if for some $p^d$-dimensional vector ${\widehat u}$
satisfying $u=P{\widehat u}$ there exists a base $B\in\B(M)$ satisfying
${\widehat W}(B)={\widehat u}$. Since ${\widehat W}$ is
$\{0,1\}$-valued, any such vector ${\widehat u}={\widehat W}(B)$
must lie in $\{0,\dots,m\}^{p^d}$, where $m:=\rank(M)$.
Therefore, checking if there is a base $B\in\B(M)$ with $W(B)=u$,
reduces to going over all vectors ${\widehat u}\in
\{0,\dots,m\}^{p^d}$, and for each, checking if  $P{\widehat u}=u$
and if there is a base $B$ satisfying ${\widehat
W}(B)={\widehat u}$. This justifies the following algorithm.

\begin{algorithm2e}[H]
{\bf let} $P$ be the $d\times p^d$ pattern matrix,
and {\bf let} $\widehat{W}$ be the selector matrix (both determined by $W$)\;
\For{$\widehat{u}\in\{0,1,\ldots,m\}^{p^d}$}{
\If {$P\widehat{u}=u$}{
\lIf {there is a base $B\in\B(M)$ with $\widehat{W}(B)=\widehat{u}$}{{\bf return} $B$\;}
}
}
{\bf return} NO;
\end{algorithm2e}

Since $d$ and $p$ are fixed and $m\leq n$, the number $(m+1)^{p^d}$ of
such potential vectors ${\widehat u}$ is polynomial in the data. For
each such vector ${\widehat u}$, checking if $P{\widehat u}=u$ is
easily done by direct multiplication; and checking if there exists a base $B\in\B(M)$
satisfying ${\widehat W}(B)={\widehat u}$ can be done in polynomial time using the
algorithm of Lemma \ref{Disjoint}
applied to the matroid $M$,
the ${\hat d}:=p^d$ weight
vectors $\widehat{w}_v$, $v\in V$\break (the $\{0,1\}$-valued rows
of the matrix $\widehat W$, having pairwise disjoint supports), and
the vector $\widehat u$. \epr

We are now in position to solve the nonlinear optimization problem
over  a matroid that is  presented by  an independence oracle.

\vskip.2cm\noindent{\bf Theorem \ref{Arbitrary}\ }
For every fixed $d$ and $p$, there is an algorithm that, given
a matroid $M$ presented by an independence oracle on  the $n$-element ground set $N$,
integers $a_1,\dots,a_p\in\Z$,
weight vectors $w_1,\dots,w_d\in\{a_1,\dots,a_p\}^n$, and function
$f:\R^d\rightarrow\R$ presented by a comparison oracle, solves the nonlinear
matroid optimization problem in time that is polynomial in $n$ and $\max \l a_i\r$.

\vskip.4cm \bpr
We have three major steps:
\begin{enumerate}

\item First, under the hypotheses of the theorem,
assume that in polynomial time we can calculate the optimal objective
value of the problem.
Then, for every subset $S\subseteq N$, the optimal
objective value of the subproblem on $S$ can be computed in
polynomial time, since the entries of the restrictions of
$w_1,\dots,w_d$ to $S$ also attain values in $\{a_1,\dots,a_p\}$.
By Lemma \ref{Objective} this implies that an optimal base can be
found and the nonlinear matroid optimization problem solved in
polynomial time. So, it remains to show that in polynomial time we can calculate
the optimal objective value of the problem.

\item
To accomplish this, we first show how to compute the set $U$ of $W$-profiles
of bases of $M$, by working with an appropriately defined superset $Z$ of $U$.
Let $m:=\rank(M)$. Consider any $i=1,\dots,d$ and
any $m$-subset $B$ of $N$. Then, since $w_{i,j}\in\{a_1,\dots,a_p\}$ for all $j$,
we have that, for some nonnegative integer vector
$\lambda_i=(\lambda_{i,1},\dots,\lambda_{i,p})$ with
$\sum_{k=1}^p\lambda_{i,k}=m$,
$$w_i(B)\ =\ \sum_{j\in B} w_{i,j}\ =\
\sum_{k=1}^p\lambda_{i,k} a_k\ =\ \lambda_i \cdot a\ ,$$
where $a:=(a_1,\dots,a_p)$. Therefore, we find that set the
set $U$ of $W$-profiles of bases satisfies
\begin{eqnarray*}
U & = & \{W(B)\ :\ B\in\B(M)\} \\
  & \subseteq & \{W(B)\ :\ B\subseteq N,\ |B|=m\} \\
   & \subseteq & Z\ :=\ \left\{(\lambda_1\cdot a ,\dots,\lambda_d\cdot a)\ :\
\lambda_i\in\{0,1,\dots,m\}^p\,,\ \ \lambda_i\cdot{\bf 1}=m\,,\ \ i=1,\dots,d\right\}\ .
\end{eqnarray*}
These observations justify the following algorithm to compute the
set $U$ of $W$-profiles of bases:

\parbox{6.25in}{
\begin{algorithm2e}[H]
Compute $m:=\rank(M)$ and {\bf let} $a:=(a_1,\dots,a_p)$\;
Start with $Z:=\emptyset$\;
\For{$\Lambda\in \{0,1,\dots,m\}^{d\times p}$}{
\lIf {$\Lambda {\bf 1}=m{\bf 1}$} {{\bf let} $Z:=Z\cup\{\Lambda a\}$\;}}
Start with $U:=\emptyset$\;
\For{$u\in Z$}{
\lIf {there is a base $B\in\B(M)$ with $W(B)=u$}{{\bf let} $U:=U\cup\{u\}$\;}
}
{\bf return} $U$\;
\end{algorithm2e}
}

Observe that, since $d$ and $p$ are fixed,
$|Z|\leq|\{0,1,\dots,m\}^{d\times p}|=(m+1)^{pd}$
is polynomially bounded in $m\le n$ and hence so are the numbers of iterations in each
of the ``for" loops of the algorithm. Also note that, in each iteration of the
second loop, we can apply the algorithm of Lemma \ref{Unary} to determine, in
polynomial time, whether there is a base
$B\in\B(M)$ with $W(B)=u$. Therefore, we can efficiently determine $U$.

\item
By repeatedly querying the comparison oracle of $f$ on $|U|-1$ suitable
pairs of points in $U$, we obtain $f^*:=\min\{f(u)\ :\ u\in U\}$, which by
Proposition \ref{BaseWeightVector} is the optimal objective value.
\end{enumerate}

Finally, by embedding steps 2-3 above as subroutines to solve the $n+1$ restrictions of
the problem to suitable subsets $S\subseteq N$ in step 1, we obtain our desired algorithm.
\epr

\section{Vectorial Matroids}

In this section we develop an algebraic algorithm for nonlinear matroid
optimization over vectorial matroids. It runs in time polynomial in the binary
length $\l A\r$ of the matrix $A$ representing the matroid and in the
unary length $\max|w_{i,j}|$ of the weights. It is much more efficient than
the combinatorial algorithm of \S3 and applies to weights with an unlimited
number of different values $w_{i,j}$ of entries.

First, we show that it suffices to deal with nonnegative weight vectors.

\bl{Nonnegative} The nonlinear matroid optimization problem with arbitrary
integer weight vectors\break $w_1\dots,w_d\in\Z^n$ is polynomial-time reducible
to the special case of nonnegative vectors $w_1\dots,w_d\in\N^n$.
\el

\bpr
Consider a matroid $M$, integer weights $w_1,\dots,w_d\in\Z^n$,
and function $f:\R^d\rightarrow\R$. Let $m:=\rank(M)$
and $q:=\max|w_{i,j}|$. Define nonnegative weights by $w'_{i,j}:=w_{i,j}+q$
for all $i,j$, and define a new function $f':\R^d\rightarrow\R$ by
$f'(y_1,\dots,y_d):=f(y_1-mq,\dots,y_d-mq)$ for every $y=(y_1,\dots,y_d)\in\R^d$.
Note that the unary length of each new weight $w'_{i,j}$ is at most
twice the maximum unary length of the original weights  $w_{i,j}$, and a
comparison oracle for $f'$ is easily realizable from a comparison oracle for $f$.
Then for every base $B$ and for each $i=1,\dots,d$, we have
$$w'_i(B)\ =\ \sum_{j\in B}w'_{i,j}\ =\ \sum_{j\in B}(w_{i,j}+q)
\ =\ \left(\sum_{j\in B}w_{i,j}\right)+mq\ =\ w_i(B)+mq\ ,$$
implying the following equality between the new and original objective function values,
$$
f'(w'_1(B),\dots,w'_d(B))\ =\ f'(w_1(B)+mq,\dots,w_d(B)+mq)\ =\ f(w_1(B),\dots,w_d(B))\ .
$$
Therefore, a base $B\in\B(M)$ is optimal for the nonlinear matroid
optimization problem with data $M$, $w_1,\dots,w_d$ and $f$ if and only if it
is optimal for the problem with $M$, $w'_1,\dots,w'_d$ and $f'$.
\epr

So we assume henceforth that the weights are nonnegative. As explained in
\S\ref{Preparation} and carried out in \S3, we will filter the set
$U=\{W(B)\ :\ B\in\B(M)\}$ of $W$-profiles of bases out of a suitable superset $Z$.
However, instead of checking if $u\in U$ for one point $u\in Z$ after the other,
we will filter here the entire set $U$ out of $Z$ at once.
We proceed to describe this procedure.

Let $A$ be a (nonempty) $m\times n$ integer matrix $A$ of full row rank $m$, and
let $M$ be the vectorial matroid of $A$. Note that $m=\rank(M)$.
Let $w_1,\dots,w_d\in\N^n$
be nonnegative integer weight vectors and let $q:=\max w_{i,j}$.
Then for each $i=1,\dots,d$ and each $m$-subset $B$ of $N$, we have
$w_i(B)\in\{0,1,\dots,mq\}$, and therefore
\begin{eqnarray*}
U & = & \{W(B)\ :\ B\in\B(M)\} \ \ \subseteq \ \ \{W(B)\ :\ B\subseteq N,\ |B|=m\}\\
   & \subseteq & Z \ \ :=\ \ \left\{0,1,\dots,mq\right\}^d \ \ \subseteq \N^d~.
\end{eqnarray*}
We will show how to filter the set $U$ out of the above superset $Z$
of potential $W$-profiles of bases. For each base $B\in\B(M)$, let
$A^B$ denote the nonsingular $m\times m$ submatrix of $A$ consisting
of those columns indexed by $B\subseteq N$. Define the
following polynomial in $d$ variables $y_1,\dots,y_d$:
\begin{equation}\label{e_polynomial}
g\quad=\quad g(y)\quad:=\quad\sum_{u\in Z}g_u y^u\quad:=\quad
\sum_{u\in Z}g_u\prod_{k=1}^d y_k^{u_k}\quad,
\end{equation}
where the coefficient $g_u$ corresponding to $u\in Z$
is the nonnegative integer
\begin{equation}\label{e_coefficient}
g_u\quad:=\quad\sum \left\{\det\!^2(A^B)\,:
\,B\in\B(M),\ \ W(B)=u\right\}\quad.
\end{equation}
Now, $\det\!^2(A^B)$ is positive for every base $B\in\B(M)$. Thus, the
coefficient $g_u$ corresponding to $u\in Z$ is nonzero if and only
if there exists a matroid base $B\in\B(M)$ with $W(B)=u$. So the
desired set $U$ is precisely the set of exponents of monomials $y^u$
with nonzero coefficient in $g$. We record this for later use:

\bp{Coefficient} Let $M$ be the vectorial matroid of an $m\times
n$ matrix $A$ of rank $m$, let $w_1,\dots,w_d\in\N^n$, and let
$g(y)$ be the polynomial in (\ref{e_polynomial}). Then
$U:=\{W(B)\ :\ B\in\B(M)\}=\{u\in Z\ :\ g_u\neq 0\}$. \ep
By Proposition
\ref{Coefficient}, to compute $U$, it suffices to compute all
coefficients $g_u$. Unfortunately, they cannot be computed directly
from the definition (\ref{e_coefficient}), since this involves again
checking exponentially many $B\in \B(M)$ -- precisely what we
are trying to avoid! Instead, we will compute the
$g_u$ by interpolation. However, in order to do so, we need a way of
evaluating $g(y)$ under numerical substitutions. We proceed to show
how this can be efficiently accomplished.

Let $Y$ be the $n\times n$ diagonal matrix whose $j$-th diagonal component
is the monomial $\prod_{i=1}^d y_i^{w_{i,j}}$ in the variables
$y_1,\dots,y_d$; that is, the matrix of monomials defined by
$$Y\quad:=\quad \diag\left(
\prod_{i=1}^d y_i^{w_{i,1}},\dots,\prod_{i=1}^d y_i^{w_{i,n}}\right)~.$$
The following lemma will enable us to compute the value of $g(y)$
under numerical substitutions.
\bl{Identity}
For any $m\times n$ matrix $A$ of rank $m$ and
nonnegative weights $w_1,\dots,w_d\in\N^n$ we have
$$g(y)\ =\ \det(AYA^T)\ .$$
\el
\bpr
By the classical Binet-Cauchy identity, for any
two $m\times n$ matrices $C,D$ of rank $m$ we have
$\det(CD^T)=\sum\{\det(C^B)\det(D^B):B\in\B(M)\}$.
Applying this to $C:=AY$ and $D:=A$, we obtain
\begin{eqnarray*}
\det(AYA^T)& =& \sum_{B\in\B(M)}\det\left((AY)^B\det(A^B)\right)\ =
\ \sum_{B\in\B(M)}\prod_{j\in B}\prod_{i=1}^d y_i^{w_{i,j}}\det(A^B)\det(A^B)\\
& = & \sum_{B\in\B(M)}\prod_{i=1}^d y_i^{w_i(B)}\det\!^2(A^B)
\ = \ \sum_{u\in Z}\sum_{B\in\B(M)\, : \atop W(B)=u}
\det\!^2(A^B)\prod_{i=1}^d y_i^{u_i}\\
& =& \sum_{u\in Z}g_u y^u\ =\ g(y)\ \ .
\end{eqnarray*}
\hfill\epr

Lemma \ref{Identity} paves the way for computing the coefficients of the polynomial
$g(y)=\sum_{u\in Z}g_u y^u$ by interpolation. We will choose sufficiently many
suitable points on the moment curve in $\R^Z$, substitute each point into $y$,
and evaluate $g(y)$ using the lemma.
We will then solve the system of linear equations for the coefficients $g_u$.
The next lemma describes the details and shows that this can be done efficiently.

\bl{Interpolation}
For every fixed $d$, there is an algorithm that, given any $m\times n$ matrix $A$
of rank $m$ and weights $w_1,\dots,w_d\in\N^n$, computes all coefficients $g_u$ of
$g(y)$ in time polynomial in $\max w_{i,j}$ and $\l A\r $.
\el
\bpr
Let $q:=\max w_{i,j}$ and $p:=mq+1$. Then a superset of potential
$W$-profiles of bases is $Z:=\{0,1,\dots,mq\}^d$ and satisfies $|Z|=p^d$.
For $t=1,2,\dots,p^d$, let $Y(t)$ be the numerical matrix
obtained from $Y$ by substituting $t^{p^{i-1}}$ for $y_i$, $i=1,\dots,d$.
By Lemma \ref{Identity} we have $g(y)=\det(AYA^T)$, and therefore we obtain
the following system of $p^d$ linear equations in the $p^d$ variables $g_u$, $u\in Z$:
\begin{eqnarray*}
\det(AY(t)A^T) &
= & \det\left( A\ \diag_j\left( \prod_{i=1}^d t^{w_{i,j}p^{i-1}} \right) A^T\right)
\quad =\quad\sum_{u\in Z}g_u\prod_{i=1}^d t^{u_i p^{i-1}}\\
&=&\sum_{u\in Z}t^{\sum_{i=1}^d u_i p^{i-1}}\ g_u~,
\quad\quad t=1,2,\dots,p^d\quad.
\end{eqnarray*}
As $u$ runs through $Z$, the sum $\sum_{i=1}^d u_i p^{i-1}$
attains precisely all $|Z|=p^d$ distinct values $0,1,\dots,p^d-1$.
This implies that, under the total order of the points $u$ in $Z$ by increasing
value of $\sum_{i=1}^d u_i p^{i-1}$, the vector of coefficients of the
$g_u$ in the equation corresponding to $t$ is precisely the point
$(t^0,t^1,\dots,t^{p^d-1})$ on the moment curve in
$\R^Z\cong\R^{p^d}$. Therefore, the equations are linearly
independent and hence the system can be uniquely solved for the $g_u$.

These observations justify the following algorithm to compute the $g_u$, $u\in Z$~:

\begin{algorithm2e}[H]
Compute $m:=\rank(A)$\;
{\bf let} $q:= \max w_{i,j}$~, and {\bf let} p:=mq+1~\;
{\bf let} $Y:=\diag_j\left(\prod_{i=1}^d y_i^{w_{i,j}}\right)$~\;
\For{$t=1,2,\ldots,p^d$}{
{\bf let} $Y(t)$ be the numerical matrix obtained by substituting $t^{p^{i-1}}$ for $y_i$, $i=1,2,\ldots,d$, in $Y$\;
Compute $\det(AY(t)A^T)$\;
}
Compute and {\bf return} the unique solution $g_u$, $u\in Z$, of the square linear system:

\qquad$\det(AY(t)A^T)=\sum_{u\in Z}t^{\sum_{i=1}^d u_i p^{i-1}} g_u~,\ \ u\in Z$~\;
\end{algorithm2e}

We now show that this system can be solved in polynomial time.
First, the number of equations and indeterminates is
$p^d=(mq+1)^d$ and hence polynomial in the data.
Second, for each $i,j=1,\dots,n$ and $t=1,2,\dots,p^d$,
it is easy to see that the $(i,j)$-th entry of $AY(t)A^T$ satisfies
$$
\left|\sum_{h=1}^n a_{i,h}\prod_{k=1}^d t^{p^{k-1}{w_{k,h}}} a_{j,h}\right|
\ \leq\ \sum_{h=1}^n |a_{i,h}a_{j,h}|p^{dp^d\max w_{k,h}}\ ,$$
implying that the binary length $\l AY(t)A^T\r$ of $AY(t)A^T$
is polynomially bounded in the data as well.

It follows that  $\det(AY(t)A^T)$ can be computed in
polynomial time by Gaussian elimination for all $t$, and the system of
equations can indeed be solved for the $g_u$ in polynomial time.
We further note that the system of equations is a Vandermonde system,
so the number of arithmetic operations needed to solve it is
just quadratic in its dimensions.
\epr

We can now efficiently solve the nonlinear optimization problem
over vectorial matroids with unary weights.

\vskip.2cm\noindent{\bf Theorem \ref{Vectorial}\ }
For every fixed $d$, there is an algorithm that, given integer $m\times n$
matrix $A$, weight vectors $w_1,\dots,w_d\in\Z^n$, and
function $f:\R^d\rightarrow\R$ presented by a comparison oracle,
solves the nonlinear optimization problem over the vectorial matroid
of $A$ in time polynomial in $\l A\r$ and $\max|w_{i,j}|$.

\vskip.4cm
\bpr
Let $M$ be the vectorial matroid of $A$. Recall that linear-algebraic
operations on $A$ can be done in polynomial time, say by Gaussian elimination.
Dropping some rows of $A$ if necessary without changing $M$, we may assume that $A$
has rank $m$. An independence oracle for $M$ is readily
realizable since $S\subseteq N$ is independent in $M$ precisely when the
columns of $A$ indexed by $S$ are linearly independent. Applying, if necessary,
the procedure of Lemma \ref{Nonnegative} and adjusting the weights while at most
doubling the unary length of the maximum weight, we may also assume that
the weights are nonnegative.

We will show how to compute the optimal objective value
$f^*:=\min\{f(W(B))\ :\ B\in\B(M)\}$ in polynomial time. This will also
imply that, for every subset $S\subseteq N$, the optimal objective
value of the subproblem on $S$ can be computed in polynomial time.
By Lemma \ref{Objective} this will show that an optimal base can be
found and the nonlinear matroid optimization problem solved in
polynomial time.

Let $q:=\max w_{i,j}$ and consider the superset $Z:=\{0,1,\dots,mq\}^d$
of potential $W$-profiles of bases and the polynomial
$g(y)=\sum_{u\in Z}g_u y^u$ as defined in (\ref{e_polynomial}) and (\ref{e_coefficient}).
By Proposition \ref{Coefficient} we have
$$U \ \ = \ \ \{W(B)\ :\ B\in\B(M)\}\ \ =\ \ \{u\in Z\ :\ g_u\neq 0\}\ \ .$$
Applying now the algorithm of Lemma \ref{Interpolation}, we can compute in
polynomial time the right-hand side and hence the left-hand side,
providing the filtration of the set $U$ of $W$-profiles of bases out of $Z$.
By repeatedly querying the comparison oracle of $f$ on suitable
pairs of points in $U$, we obtain $f^*:=\min\{f(u)\ :\ u\in U\}$ which by
Proposition \ref{BaseWeightVector} is the desired optimal objective value.
\epr

\section{Experimental Design}

We now discuss applications of nonlinear matroid optimization to experimental design.
For general information on experimental design see e.g. the monograph \cite{PRW}
and references therein. As outlined in the introduction, we consider the following
rather general framework. We wish to learn an unknown system whose output
$y$ is an unknown function $\Phi$ of a multivariate input $x=(x_1,\dots,x_k)\in\R^k$.
It is customary to call the input variables $x_i$ {\em factors} of the system.
We perform several experiments. Each experiment $i$
is determined by a point $p_i=(p_{i,1},\dots,p_{i,k})$ and consists of feeding
the system with input $x:=p_i\in\R^k$ and measuring the corresponding output
$y_i:=\Phi(p_i)\in\R$. Based on these experiments, we wish to {\em fit a model}
for the system, namely, determine an estimation $\hat{\Phi}$ of $\Phi$, that:
\begin{itemize}
\item
Lies in a prescribed class of functions;
\item
Is consistent with the outcomes of the experiments;
\item
Minimizes the {\em aberration} - a suitable criterion -
among models in the class.
\end{itemize}

We concentrate on (multivariate) {\em polynomial models} defined as follows.
Each nonnegative integer vector $\alpha\in\N^k$ serves as an
{\em exponent} of a corresponding monomial $x^\alpha:=\prod_{h=1}^k x_h^{\alpha_h}$
in the system input $x\in\R^k$. Each finite subset $B\subset \N^k$
of exponents provides a {\em model} for the system, namely a polynomial
{\em supported on $B$}, i.e. having monomials with exponents in $B$ only,
\[
\Phi_B(x)\ =\ \sum_{\alpha\in B} c_{\alpha} x^{\alpha},
\]
where the $c_{\alpha}$ are real coefficients that
need to be determined from the measurements by interpolation.

We assume that the set of design points $\{p_1,\dots,p_m\}\subset \R^k$ is prescribed.
Indeed, in practical applications, it may be impossible or too costly to conduct
experiments involving arbitrarily chosen points. The problem of choosing the design
(termed the {\em inverse problem} in the statistics literature, see \cite{ADT}
and references therein), is of interest in its own right,
and its computational aspects will be considered elsewhere.
We collect the design points in an $m\times k$ \emph{design} matrix $P$.
Thus, the $i$-th row of this matrix is the $i$-th design point $p_i$.
A model $B\subset\N^k$ is {\em identifiable} by a design $P$
if for any possible measurement values $y_i=\Phi(p_i)$ at the design points,
there is a unique polynomial $\Phi_B(x)$ supported on $B$ that interpolates $\Phi$,
that is, satisfies $\Phi_B(p_i)=y_i=\Phi(p_i)$ for every
design point $p_i=(p_{i,1},\dots,p_{i,k})$.

Among models identifiable by a given design, we wish to determine one that is best
under a suitable criterion. Roughly speaking, common criteria
ask for {\em low degree polynomials}. To make this precise, for each identifiable
model $B$, consider the following {\em average degree vector} whose $i$-th entry is
the average degree of variable $x_i$ over all monomials supported on $B$~:
\[
{1\over|B|}\sum_{\alpha\in B}\alpha\ =\ \left({1\over|B|}\sum_{\alpha\in B}\alpha_1,
\dots,{1\over|B|}\sum_{\alpha\in B}\alpha_k\right)~.
\]
Now, given any function $f:\R^k\rightarrow\R$,
the {\em aberration} of model $B$ induced by $f$ is defined to be
\[
{\mathcal A}(B)\ :=\ f\left({1\over|B|}\sum_{\alpha\in B}\alpha\right)~.
\]
The term {\em aberration} is the one used in the statistics literature
in this context, see e.g. \cite{FH,WW} and the reference therein.
We now give some concrete examples of functions providing useful aberrations.

\be{Aberration1}{\bf Some concrete useful aberrations.}
\begin{itemize}

\item
Consider the function $f(x):=x_1+\cdots+x_k$.
Then the aberration of model $B$ is
\[
{\mathcal A}(B)\ = \ {1\over |B|}\sum_{\alpha\in B}\alpha_1+
\cdots +{1\over |B|}\sum_{\alpha\in B}\alpha_k
\ = \ {1\over |B|}\sum_{\alpha\in B} \sum_{i=1}^k\alpha_i~,
\]
which is the average total degree of monomials supported on $B$.

\item
Consider $f(x):=\pi_1x_1+\cdots+\pi_kx_k$
for some real weights $\pi_1,\dots,\pi_k$. Then
$${\mathcal A}(B)\ = \ \pi_1{1\over |B|}\sum_{\alpha\in B}\alpha_1+
\cdots +\pi_k{1\over |B|}\sum_{\alpha\in B}\alpha_k$$
is the weighted average degree, allowing for preferences
of some variables over others.

\item
Consider the function $f(x):=\max\{x_1,\dots,x_k\}$.
Then the aberration of $B$ is
$${\mathcal A}(B)\ = \ \max\left\{{1\over |B|}\sum_{\alpha\in B}\alpha_1,
\dots,{1\over |B|}\sum_{\alpha\in B}\alpha_k\right\}\ ,$$
and is the maximum over variables of the average
variable degree of monomials supported on $B$.

\item
More generally, consider
$f(x):=\|\pi\cdot x\|_q=(\sum_{i=1}^k |\pi_i x_i|^q)^{1\over q}$.
Then the aberration of $B$ is
$${\mathcal A}(B)\ =\ \left\|\pi\cdot
\left({1\over |B|}\sum_{\alpha\in B} \alpha \right) \right\|_q\ ,$$
which is the $l_q$-norm of the weighted average degree
vector of monomials supported on $B$.
\end{itemize}
\ee

\vskip.2cm
We can now formally define the minimum-aberration model-fitting problem.

\vskip.2cm\noindent
{\bf Minimum-Aberration Model-Fitting Problem}. Given a design
$P=\{p_1,\dots,p_m\}$ of $m$ points in $\R^k$, a set
$N=\{\beta_1,\dots,\beta_n\}$ of $n$ potential exponents in $\N^k$,
and a function $f:\R^k\rightarrow\R$, find a model $B\subseteq N$
that is identifiable by $P$ and is of minimum aberration
$${\mathcal A}(B)\ =\ f\left({1\over{|B|}}\sum_{\beta_j\in B}\beta_j\right)~.$$

\vskip.2cm
It can be verified that a model $B$ is identifiable by a design $P$ if and only if
$B=\{\beta_1,\dots,\beta_m\}\subset\N^k$ for some $m$ exponents
$\beta_j=(\beta_{j,1},\dots,\beta_{j,k})$, $j=1,\dots,m$,
and the $m\times m$ matrix $T$, defined by
\[
t_{i,j}\ :=\ p_i^{\beta_j}\ =\ \prod_{h=1}^k p_{i,h}^{\beta_{j,h}}~,
\]
is invertible. (In the terminology of algebraic geometry,
the model $B$ is identifiable if the congruence classes of the monomials
$x^{\beta_1},\dots,x^{\beta_m}$ form a basis for the quotient of
the algebra of polynomials $\R[x_1,\dots,x_k]$ modulo the ideal of
polynomials vanishing on the design points; see \cite{BOT,OS} and
references therein for more on this.) If $B$ is identifiable then,
given any vector of measurements $y\in\R^m$ at the design
points, the vector of coefficients $c\in\R^m$ of the unique
polynomial $\Phi_B(x)=\sum_{j=1}^m c_j \prod_{h=1}^k x_h^{\beta_{j,h}}$
supported on $B$ that is consistent with the measurements, is given by $c:=T^{-1}y$.

\vskip.2cm
It should be now quite clear how to formulate the minimum-aberration
model-fitting problem over a design $P=\{p_1,\dots,p_m\}$, a set
$N=\{\beta_1,\dots,\beta_n\}$ of potential exponents, and a function
$f:\R^k\rightarrow\R$, as a nonlinear matroid optimization problem.
Consider the following $m\times n$ matrix $A$ having rows indexed by $P$
and columns indexed by $N$, defined by
\[
a_{i,j}\ :=\ p_i^{\beta_j}\ =\ \prod_{h=1}^k p_{i,h}^{\beta_{j,h}}~,
\quad i=1,\dots,m\ ,\ j=1,\dots,n~.
\]
Then a model $B\subseteq N$ is identifiable by $P$ if and only if the submatrix
of $A$ comprising columns indexed by $B$ is invertible. If the rank of $A$
is less than $m$, then no $B\subseteq N$ is identifiable, and the set $N$ of
potential exponents should be augmented with more exponents.
So assume that $A$ has rank $m$. Let $M$ be the vectorial matroid
of $A$, so that
\[
\B(M)\ :=\ \{B\subseteq N\ :\  B\ \mbox{is identifiable by}\ P\}~.
\]
Now define $k$ weights vectors $w_1,\dots,w_k\in\N^n$ by
$w_{i,j}:={1\over m}\beta_{j,i}$ for $i=1,\dots,k$, $j=1,\dots,n$.
Then the aberration of model $B$ is
$${\mathcal A}(B)\ = \
f\left({1\over m}\sum_{\beta_j\in B}\beta_j\right)\ =\
f\left(\sum_{\beta_j\in B} w_{1,j},\dots,\sum_{\beta_j\in B} w_{k,j}\right)\ .$$
Thus, the aberration of a model $B$
identifiable by the design $P$ is precisely the objective function
value of the base $B$ in the nonlinear matroid optimization problem
over the matroid $M$ above, with $d:=k$ and the weights
$w_1\dots,w_k\in\N^n$ as above, and with the given function
$f:\R^k\rightarrow \R$. Assuming that a comparison oracle for
the function $f$ can be realized, which is practically always true,
and that the design points are rational so that they can be input
and processed on a digital computer, we obtain the following
corollary of Theorems \ref{Arbitrary} and \ref{Vectorial}:

\vskip.2cm\noindent{\bf Corollary \ref{ExperimentalDesign}\ }
For every fixed $k$, there is an algorithm that, given a rational design
$P=\{p_1,\dots,p_m\}$ in $\R^k$, a set $N=\{\beta_1,\dots,\beta_n\}$
in $\N^k$, and a function $f:\R^k\rightarrow\R$ presented by a
comparison oracle, solves the minimum-aberration model-fitting problem
in time polynomial in $m,n$, $\l P\r$ and $\max \beta_{i,j}$.

\vskip.2cm
It is very natural and common in practice to consider {\em hierarchical models};
that is, models $B$ with the property that $\beta\leq\alpha\in B$ implies $\beta\in B$.
In \cite{BOT,OS}, in the context of the theory of Gr\"obner bases in commutative algebra,
it was shown that the smallest set containing all $m$-point hierarchical models $B$
(termed {\em staircases} therein) in $\N^k$ is the following set, consisting
of roughly $O(m\log m)$ points,
$$N\ :=\ \{\alpha\in\N^k\ :\ \prod_{h=1}^k(\alpha_h+1)\leq m\}\ .$$
Thus, Corollary \ref{ExperimentalDesign} will be typically applied
with this set $N$ as the set of potential monomial exponents.

\vskip.5cm
We proceed to describe a more general useful class of aberrations that can be treated,
which is naturally suggested by the nonlinear matroid optimization formulation.
As before, we are given a design $P=\{p_1,\dots,p_m\}$ in $\R^k$ and a set
$N=\{\beta_1,\dots,\beta_n\}$ of potential exponents in $\N^k$.
But now we are also given $d$ weight vectors $w_1,\dots,w_d\in\Z^n$.
The function $f$ is now defined on $\R^d$ rather than $\R^k$. The aberration
induced by the weights and the function is now simply the objective function
of the nonlinear matroid optimization problem, which for an identifiable
model $B\subseteq N$ is given by
$${\mathcal A}(B)\ :=\
f\left(W(B)\right)\ ,$$
where $W$ is the matrix with rows $w_i$. Note that aberrations of the type
considered before can be recovered as a special case with $d:=k$ and
$w_{i,j}:={1\over m}\beta_{j,i}$ for all $i,j$. Here are a few useful examples.

\be{Aberration2}{\bf Some concrete useful generalized aberrations.}
Let $N=\{\beta_1,\dots,\beta_n\}\subset\N^k$ be any set of exponents.
Let $\theta$ be a small positive integer, say $\theta=1$ or $\theta=2$,
that will serve as a desired bound on the degrees of variables in monomials
of the sought after models $B$ contained in $N$.
\begin{itemize}

\item
Let $d:=1$, and define the single weight vector $w_1\in\N^n$ by
$$w_{1,j}\quad:=\quad\left\{
\begin{array}{ll}
    0, & \hbox{if $\beta_{i,j}\leq \theta$~, for every $i=1,\dots,k$~;} \\
    1, & \hbox{otherwise,} \\
\end{array}
\right.
$$
for $j=1,\dots,n$.
Then $\sum_{\beta_j\in B}w_{1,j}$ is the number of monomials supported on
$B\subseteq N$ that do {\em not} meet the degree bound; in particular, for $\theta=1$
it is the number of non square-free monomials. Taking $f:\R\rightarrow\R$
to be the identity $f(y):=y$, the aberration ${\mathcal A}(B)$ of model $B$
is the number of undesired monomials. In particular, an optimal model $B$
has ${\mathcal A}(B)=0$ if and only if the design admits an identifiable
model with all variables in all monomials having degree at most $\theta$.

\item
Now let $d:=k$, and define weight vectors $w_1,\dots,w_k\in\N^n$ by
$$w_{i,j}\quad:=\quad\left\{
\begin{array}{ll}
    0, & \hbox{if $\beta_{i,j}\leq \theta$~;} \\
    1, & \hbox{otherwise,} \\
\end{array}
\right.
$$
for $i=1,\dots,k$, $j=1,\dots,n$.
Then $\sum_{\beta_j\in B}w_{i,j}$ is the number of monomials supported
on model $B\subseteq N$ for which variable $x_i$ violates the degree bound $\theta$.
Defining $f:\R^k\rightarrow\R$ by $f(y):=\max_{i=1}^k y_i$,
we get that the aberration ${\mathcal A}(B)$ of model $B$ is the maximum over
variables of the number of monomials having $x_i$ violating the degree
bound $\theta$. The optimal model will minimize the maximum violation.
\end{itemize}
\ee

We have the following generalized minimum-aberration model-fitting problem and corollary.

\vskip.2cm\noindent
{\bf Generalized Minimal-Aberration Model-Fitting Problem}.
Given a design $P=\{p_1,\dots,p_m\}$ in $\R^k$, a set $N=\{\beta_1,\dots,\beta_n\}$
of potential exponents in $\N^k$, weight vectors $w_1,\dots,w_d\in\Z^n$,
and a function $f:\R^d\rightarrow\R$,
find a model $B\subseteq N$ that is identifiable by $P$ and is of minimum aberration
$${\mathcal A}(B)\ =\ f\left(W(B)\right)\ =\
f\left(\sum_{\beta_j\in B} w_{1,j},\dots,\sum_{\beta_j\in B} w_{d,j}\right)\ .$$

\bc{GeneralizedExperimentalDesign}
For every fixed $k$ and $d$, there is an algorithm that, given any rational design\break
$P=\{p_1,\dots,p_m\}$ in $\R^k$, any set $N=\{\beta_1,\dots,\beta_n\}$ of potential
exponents in $\N^k$, weight vectors $w_1,\dots,w_d\in\Z^n$, and function
$f:\R^d\rightarrow\R$ presented by a comparison oracle, solves the generalized
minimum-aberration model-fitting problem in time polynomial in $m,n$, $\l P\r$,
$\max \beta_{j,i}$ and $\max |w_{i,j}|$.
\ec

\section*{Acknowledgements}

The research of Shmuel Onn and Henry Wynn was partially supported
by the Joan and Reginald Coleman-Cohen Exchange Program during a
stay of Henry Wynn at the Technion. The research of Jon Lee,
Shmuel Onn and Robert Weismantel was partially supported by the
Mathematisches Forschungsinstitut Oberwolfach during a stay
within the Research in Pairs Programme.
Yael Berstein was supported by an Irwin and Joan Jacobs Scholarship
and by a scholarship from the Graduate School of the Technion.
Hugo Maruri-Aguilar and Henry Wynn were also supported by the Research Councils
UK (RCUK) Basic Technology grant ``Managing Uncertainty in Complex Models".
Shmuel Onn was also supported by the ISF - Israel Science Foundation.
Robert Weismantel was also supported by the European TMR Network ADONET 504438.

\noindent {\small Yael Berstein}\newline
\emph{Technion - Israel Institute of Technology, 32000 Haifa, Israel}\newline
\emph{email: yaelber{\small @}tx.technion.ac.il}

\vskip.3cm\noindent {\small Jon Lee}\newline
\emph{IBM T.J. Watson Research Center, Yorktown Heights, NY 10598, USA}\newline
\emph{email: jonlee{\small @}us.ibm.com},
\ \ \emph{http://www.research.ibm.com/people/j/jonlee}

\vskip.3cm\noindent {\small Hugo Maruri-Aguilar}\newline
\emph{London School of Economics, London WC2A 2AE, UK}\newline
\emph{email: h.maruri-aguilar{\small @}lse.ac.uk}

\vskip.3cm\noindent {\small Shmuel Onn}\newline
\emph{Technion - Israel Institute of Technology, 32000 Haifa, Israel}\newline
\emph{email: onn{\small @}ie.technion.ac.il},
\ \ \emph{http://ie.technion.ac.il/{\small $\sim$onn}}

\vskip.3cm\noindent {\small Eva Riccomagno}\newline
\emph{Universita' degli Studi di Genova, 16146 Genova, Italia}\newline
\emph{email: riccomagno{\small @}dima.unige.it},
\ \ \emph{http://www.dima.unige.it/{\small $\sim$riccomagno}}

\vskip.3cm\noindent {\small Robert Weismantel}\newline
\emph{Otto-von-Guericke Universit\"at Magdeburg,
D-39106 Magdeburg, Germany}\newline
\emph{email: weismantel{\small @}imo.math.uni-magdeburg.de},
\ \ \emph{http://www.math.uni-magdeburg.de/{\small $\sim$weismant}}

\vskip.3cm\noindent {\small Henry Wynn}\newline
\emph{London School of Economics, London WC2A 2AE, UK}\newline
\emph{email: h.wynn{\small @}lse.ac.uk}

\end{document}